\documentclass[10]{article}
\usepackage{fleqn}
\usepackage{graphicx}
\usepackage{amsfonts,amsthm}
\usepackage{amsthm} 
\usepackage{amsmath}
\usepackage{amssymb}
\usepackage{AlgorithmDefinition}
\usepackage{lstalgorithms}
\usepackage{hyperref}
\def\PP{\mathbb{P}} 
\def\PG{{\rm PG}}
\theoremstyle{plain}

\theoremstyle{definition}

\theoremstyle{remark}

\begin{document}
\Large
\begin{center}
{\bf Veldkamp Spaces of Low-Dimensional Ternary\\ Segre Varieties}
\end{center}
\vspace*{-.1cm}
\large
\begin{center}
 Metod Saniga,$^{1}$ J\'er\^ome Boulmier,$^{2}$ Maxime Pinard$^{2}$ and Fr\'ed\'eric Holweck$^{3}$
\end{center}
\vspace*{-.5cm} \normalsize
\begin{center}

$^{1}$Astronomical Institute of the Slovak Academy of Sciences,\\
SK-05960 Tatransk\' a Lomnica, Slovak Republic\\
(msaniga@astro.sk)

\vspace*{.2cm}

$^{2}$Universit\'e de Technologie de Belfort-Montb\'eliard,
F-90010 Belfort, France\\ (jerome.boulmier@utbm.fr, maxime.pinard@utbm.fr)
%\vspace*{0.4cm}

and

$^{3}$ICB/UTBM, University of Bourgogne-Franche-Comt\'e
F-90010 Belfort, France\\(frederic.holweck@utbm.fr)

\end{center}

\vspace*{-.3cm} \noindent \hrulefill

\vspace*{-.0cm} \noindent {\bf Abstract}

\noindent Making use of the `Veldkamp blow-up' recipe, introduced by Saniga and others (Ann. Inst. H. Poincar\' e D2 (2015) 309) for binary Segre varieties, we study geometric hyperplanes and Veldkamp lines of Segre varieties $S_k(3)$, where  $S_k(3)$  stands for the $k$-fold direct product of projective lines of size four and $k$ runs from 2 to 4. Unlike the binary case, the Veldkamp spaces here feature also non-projective elements. Although for $k=2$ such elements are found only among Veldkamp lines, for $k \geq 3$ they are also present among Veldkamp points of the associated Segre variety. Even if we consider only projective geometric hyperplanes, we find four different types of non-projective Veldkamp lines of $S_3(3)$, having 2268 members in total, and five more types if non-projective ovoids are also taken into account.
Sole geometric and combinatorial arguments lead to as many as 62 types of projective Veldkamp lines of $S_3(3)$, whose blowing-ups yield 43 distinct types of projective geometric hyperplanes of $S_4(3)$. As the latter number falls short of 48, the number of different large orbits of $2 \times 2 \times 2 \times 2$ arrays over the three-element field found by Bremner and Stavrou (Lin. Multilin. Algebra 61 (2013) 986), there are five (explicitly indicated) hyperplane types such that each is the merger of two different large orbits. Furthermore, we single out those 22  types of geometric hyperplanes of $S_4(3)$, featuring 7\,176\,640 members in total, that are in a one-to-one correspondence with the points lying on the unique hyperbolic quadric $\mathcal{Q}_0^{+}(15,3) \subset {\rm PG}(15,3) \subset \mathcal{V}(S_4(3))$; and, out of them, seven ones that correspond bijectively to the set of 91\,840 generators of the symplectic polar space $\mathcal{W}(7,3) \subset \mathcal{V}(S_3(3))$. 
For $k=3$ we also briefly discuss embedding of the binary Veldkamp space into the ternary one. Interestingly, only 15 (out of 41) types of lines of 
$\mathcal{V}(S_3(2))$ are embeddable and one of them, surprisingly, into a {\it non}-projective line of $\mathcal{V}(S_3(3))$ only.

\vspace*{.3cm}

%{\bf MSC Codes:} 51Exx, 81R99\\
%{\bf PACS Numbers:} 02.10.Ox, 02.40.Dr, 03.65.Ca\\

\vspace*{-.2cm} \noindent \hrulefill

%\vspace*{-.3cm} 
\section{Introduction} 
%\vspace*{-.3cm} 
To any point-line incidence geometry $\mathcal{C}$ endowed with geometric hyperplanes one can associate a particular space $\mathcal{V}(\mathcal{C})$, called the Veldkamp space of $\mathcal{C}$, as follows \cite{buec}: a point of  $\mathcal{V}(\mathcal{C})$ is a geometric hyperplane of $\mathcal{C}$ and its line is the collection $H'H''$ of all geometric hyperplanes $H$ of $\mathcal{C}$  such that $H' \cap H'' = H' \cap H = H'' \cap H$ or $H = H', H''$, where $H'$ and $H''$ are distinct. Clearly, the structure of $\mathcal{V}(\mathcal{C})$ depends strongly on properties of $\mathcal{C}$. When $\mathcal{C}$ is a partial gamma space with all lines of size three, then \cite{coh} the lines of $\mathcal{V}(\mathcal{C})$ are of the form $(H', H'', H' \oplus H'')$, where $H' \oplus H'' \equiv (H'\cap H'') \cup (\mathcal{C} \backslash H'\cup H'')$, and
$\mathcal{V}(\mathcal{C})$ is a projective space (over the smallest Galois field). $\mathcal{V}(\mathcal{C})$ has also interesting properties when $\mathcal{C}$ exhibits a highly-symmetric layered structure \cite{near}. However, the importance of this concept can properly be appreciated when $\mathcal{C}$ features both layering and nesting of some basic structures. One of the simplest and most illustrative examples of such a geometry is the Segre variety that is a $k$-fold product of projective lines of size three, $S_k(2)$. Given the fact that  $S_k(2) =  {\rm PG}(1,2) \times S_{k-1}(2)$, where PG$(1,2)$ represents the projective line having three points, it turns out that the information about {\it points} of $\mathcal{V}(S_k(2))$, $k \geq 2$, is fully encoded in the properties of {\it lines} of $\mathcal{V}(S_{k-1}(2))$. 
The way how this encoding works, as well as its further implications, were described in detail in \cite{seg2}, which is also our basic reference. In our present paper we shall use the strategy of \cite{seg2} and perform an analogous analysis
on Segre varieties whose lines feature four points. As the complexity of this case is considerably higher when compared with the binary one, our work substantially relies on computer calculations.

The paper is organized as follows. After highlighting the relevant concepts and notation (Sec.\,2) the blow-up procedure and involved computer codes are described (Sec.\,3). The central part of the paper is Sec.\,4, which provides a detailed account of the properties of the Veldkamp spaces of $S_k(3)$, $2 \leq k \leq 4$. Finally, Sec.\,5 offers a brief summary of the main results and a passing comparison with the binary case.

%\vspace*{-.3cm} 
\section{Basic concepts and notation}
%\vspace*{-.3cm} 
Let us start with highlighting of relevant geometrical and combinatorial concepts and setting up the notation.

Given a projective space PG$(d,q)$, where $d$ is a positive integer and $q$ a power of a prime, a finite Segre variety is defined as \cite{ht}: $S_{d_1, d_2,\ldots,d_k}(q) \equiv {\rm PG}(d_1,q) \times {\rm PG}(d_2,q) \times \ldots \times {\rm PG}(d_k,q)$. Here we will only be concerned with the case $S_{1, 1,\ldots,1}(3) \equiv S_{k}(3)$ that, as an abstract point-line incidence structure, has $4^k$ points and $k 4^{k-1}$ lines, with four points per line and $k$ lines through a point. Next it comes a {\it point-line incidence structure} $\mathcal{C} = (\mathcal{P},\mathcal{L},I)$ where $\mathcal{P}$ and $\mathcal{L}$ are, respectively, sets of points and lines and where incidence $I \subseteq \mathcal{P} \times \mathcal{L}$ is a binary relation indicating which point-line pairs are incident (see, e.\,g., \cite{shult}).    A {\it geometric hyperplane} of $\mathcal{C} = (\mathcal{P},\mathcal{L},I)$ is a proper subset of $\mathcal{P}$ such that a line from $\mathcal{L}$ either lies fully in the subset, or shares with it only one point. Given a hyperplane $H$ of $\mathcal{C}$, one defines the
{\it order} of any of its points as the number of lines through the point that are fully contained in $H$; a point of maximal order  is called {\it deep}. 
Any $S_{k}(q)$, $k \geq 2$, exhibits two particular types of a geometric hyperplane: a {\it singular} hyperplane (also called a perp-set for $k=2$), comprising all the points not at maximum distance from a given point (to be called its {\it nucleus}), and an {\it ovoid}, i.\,e. a set of points partitioning the set of lines. 
As we will see later, alongside the {\it ordinary} lines of $\mathcal{V}(\mathcal{C})$ defined in Sec.\,1, for $\mathcal{C} = S_k(3)$ it will also be necessary to take into account their trivial (or extraordinary using the language of \cite{seg2}) counterparts; a {\it trivial} Veldkamp line of $S_k(3)$ consists of the   $S_k(3)$ and any of its geometric hyperplanes counted three times.
In either case, the elements shared by all the geometric hyperplanes constituting a given Veldkamp line will be called the {\it core} of this line. 

The Segre variety $S_k(3)$ (and similarly its generalization  $S_{d_1, d_2,\ldots,d_k}(q)$) can be embedded via the Segre embeding as a subgeometry of $\PG(2^k-1,3)$. Let us denote by $V^2_3$ the two-dimensional vector space over the three-element field $GF(3)$. Then  $\PG(2^k-1,3) \sim \PP(\underbrace{V_3^2\otimes \dots \otimes V_3^2}_{k \text{ times}})$ and we can define the following map 
\begin{equation}
\begin{array}{ccc}
 S_k(3)=\PG(1,3)\times\dots\times \PG(1,3) & \to & \PG(2^{k-1}-1,3)=\PP(\underbrace{V_3^2\otimes \dots \otimes V_3^2}_{k \text{ times}}),\\
 ([x_1],\dots,[x_k]) & \mapsto & [x_1\otimes\dots \otimes x_k],
\end{array}
\end{equation}
where $x_i$ denotes a vector of $V_3^2$ and $[x_i]$ the corresponding point in the projective space.

Let us assume that PG$(d,q)$ is endowed with a non-degenerate reflexive sesquilinear form or a non-singular quadratic form. Then the polar space in this ambient  PG$(d,q)$ consists of those projective subspaces that are totally isotropic or totally singular with respect to such form, respectively; projective subspaces of maximal dimensions being called generators. There exist five distinct types of finite classical polar spaces (see, e.\,g., \cite{bkm}), but we will only be concerned with (a particular type of) the {\it hyperbolic} quadric  $\mathcal{Q}_0^{+}(2n-1,q)$, $n \geq 1$, formed by all points of PG$(2n - 1, q)$ satisfying the standard equation $x_1x_2+x_3x_4+\cdots+x_{2n-1}x_{2n} = 0$, and the {\it symplectic}  polar space $\mathcal{W}(2n-1,q)$, $n \geq 1$, that consists of all points of PG$(2n-1,q)$ and all totally isotropic subspaces with respect to the standard symplectic form
$x_1y_2-x_2y_1+\cdots+x_{2n-1}y_{2n} - x_{2n}y_{2n-1}$;  $\mathcal{Q}_0^{+}(2n-1,q)$ contains $(q^{n-1}+1)(q^{n}-1)/(q-1)$ points and $\mathcal{W}(2n-1,q)$ features $(q + 1)(q^2 + 1)\cdots(q^n+1)$ generators.

Finally, a few additional remarks concerning $S_k(3)$. It possesses $4k$ $S_{k-1}(3)$'s arranged into $k$ distinct quadruples of pairwise disjoint members. Additionally,  it contains $k$ {\it distinguished} spreads of lines; a distinguished spread of lines is a set of $4^{k-1}$ mutually skew lines, each being incident with all the four $S_{k-1}(3)$'s in some quadruple. 
Also, when we will be speaking of distance between two points of $S_k(3)$ we will always have in mind the distance between the corresponding vertices in its collinearity graph.
The geometry $S_k(3)=\PG(1,3)\times \dots \times\PG(1,3)$ is stabilized by the action of the group $G=(\sigma_4\times\sigma_4\times\dots\times \sigma_4)\rtimes \sigma_k$ where $\sigma_n$ stands for the symmetric group on $n$-elements. The group $\sigma_4\times \sigma_4\times \dots \times \sigma_4$ acts naturally on $S_k(3)$, i.\,e. each $\sigma_4$ permutes points of the corresponding PG$(1,3)$, whereas  $\sigma_k$ exchanges the individual PG$(1,3)$'s.

%The elements of $\sigma_4\times \dots\sigma_4$ act by permuting the points on each projective lines of $S_3(4)$ while the action of $\sigma_k$ permute the lines of $S_k(3)$.

%\vspace*{-.5cm}
%\newpage
\section{Blow-up procedure and its computer encoding}

In this section we describe, in terms of pseudo-codes, the main algorithms that implemented, adjusted and generalized the blow-up procedure of  \cite{seg2} to generate geometric hyperplanes of  $S_k(3)$ from the knowledge of Veldkamp lines of $S_{k-1}(3)$. The codes, written in C++, are freely available at 
%\newline
\url{https://github.com/Lomadriel/HyperplaneFinder}.

We shall start by recalling the main idea that leads in \cite{seg2} to the blow-up construction. Let us consider a point-line incidence structure $X$  over the finite field $GF(q)$ such that $X$ possesses geometric hyperplanes and Veldkamp lines comprising $q+1$ points. 
Then the new incidence geometry $X\times \PG(1,q)$ also contains geometric hyperplanes, which can be obtained by 
blowing up Veldkamp lines of $X$. More precisely, let us denote by $h_0,h_1,\dots,h_q$ a set of hyperplanes of $X$ such that $\{h_0, h_1,\dots,h_q\}$ is a line of $\mathcal{V}(X)$. Then, if we denote by $[x_0],\dots,[x_q]$ the $q+1$ points of $\PG(1,q)$, one obtains a hyperplane of $X\times \PG(1,q)$ of the following type \begin{equation}\label{type1} H=\cup_i h_i\times\{[x_i]\} \subset X\times \PG(1,q).\end{equation} Permuting the points $\{[x_i]\}$ of $\PG(1,q)$ leads to $(q+1)!$ distinct geometric hyperplanes from a single Veldkamp line of $X$. %There is another type of hyperplanes in $X\times \PG(1,q)$ obtained from one single hyperplane of $X$. 
 Another type of geometric hyperplane of $X\times \PG(1,q)$ is of the form \begin{equation}\label{type2} H'=X\times \{[x_j]\} \cup_{0\leq i\neq j\leq q} h\times \{[x_i]\}\subset X\times \PG(1,q)\end{equation}
where $h$ is any hyperplane of $X$. 
As $j\in \{0,\dots,q\}$, one can generate $q+1$ different geometric hyperplanes similar to $H'$. In fact, all hyperplanes of $X\times \PG(1,q)$ can be obtained either from Eq. (\ref{type1}) or Eq. (\ref{type2}) (see \cite{seg2}). 

From the above considerations it follows that in order to find all hyperplanes of $S_k(3)$, one first needs to generate all lines\footnote{Lines producing geometric hyperplanes given by Eq. (\ref{type1}) will be called ordinary lines of $\mathcal{V}(S_{k-1}(3))$, while those corresponding to Eq. (\ref{type2}) will be called trivial ones.} of $\mathcal{V}(S_{k-1}(3))$.  This is done by Algorithm \ref{lst:algo1}, which performs the following specific tasks:
\begin{enumerate}
 \item for all pairs of distinct geometric hyperplanes $(h_i,h_j) \in (\mathcal{V}(S_{k-1}(3))^2$ find all geometric hyperplanes $h$ such that $h\cap h_i=h\cap h_j=h_i\cap h_j$ and make a list $L$ of them;
 \item for all pairs of $(h_a,h_b)\in L^2$ check if $\{h_1,h_2,h_a,h_b\}$ is a Veldkamp line, i.\,e. if $h_a\cap h_b=h_1\cap h_2$;
 \item separate the lines that are projective from those that are not.
\end{enumerate}

This last step of the algorithm is specific for the $GF(3)$-setting. In \cite{seg2} it was shown that, in light of \cite{coh}, $\mathcal{V}(S_k(2))=\PG(2^k-1,2)$, i.\,e. all lines of $\mathcal{V}(S_k(2))$ are projective lines. However, over $GF(3)$, as we will see in Sec. \ref{sec4}, through two distinct points  $h_1, h_2\in \mathcal{V}(S_{k}(3))$  there can pass more than one line. To discriminate a projective line from a non-projective one, in our computer search we applied the following test. Let $l=\{h_1,h_2,h_3,h_4\}$ be a line of $\mathcal{V}(S_{k-1}(3))$ that generates a hyperplane $H_l$ of $S_k(3)$ of type given by Eq. (\ref{type1}). Embed $S_k(3)$ into $\PG(2^{k}-1,2)$ via the Segre embedding and check if the linear span of the points of this $H_l$ in $\PG(2^k-1,2)$ fills the whole space. If so, then $H_l$ cannot be obtained as a linear section of $S_k(3)$ embedded in $\PG(2^k-1,3)$ and, hence,  $l$ is {\it not} a projective line. On the other hand, if 
the span of the embedding of $H_l$  defines a codimension-$1$ projective space in $\PG(2^k-1,3)$, i.\,e. if $H_l$ is a linear section of the Segre embedding of $S_k(3)$  in $\PG(2^k-1,3)$, then the line $l$ is projective. The test of the dimension of the span is done by computing the rank of the matrix corresponding to the span. The function ``\text{compute\_associated\_matrix}" generates the matrix of the span of the vectors in  $V_3^2\otimes V_3^2\otimes\dots\otimes V_3^2$ associated to the points of $H_l$.

\lstinputalgo[%
	label={lst:algo1},%
	caption={Veldkamp lines generation}%
]{compute_veldkamp_lines.algo}

Algorithm \ref{lst:algo2} generates from  Veldkamp lines of $S_{k-1}(3)$ geometric hyperplanes of $S_k(3)$ in terms of Eq. (\ref{type1}) and Eq. (\ref{type2}).

\lstinputalgo[%
	label={lst:algo2},%
	caption={Hyperplanes generation from Veldkamp lines of the lower dimension}%
]{compute_hyperplanes_from_veldkamp_lines.algo}

Once the geometric hyperplanes (projective or not) of $S_k(3)$ are generated, they are first classified, following \cite{seg2}, by making sole use of combinatorial and geometric arguments and their properties are displayed in form of tables. This geometro-combinatorial classification  is $(\sigma_4\times\sigma_4\times\dots\times \sigma_4)\rtimes \sigma_k$-invariant by construction.
It, however, turns out that -- as we will explicitly see in the next section for the case $k=4$ -- this is generally not enough to distinguish {\it all} the orbits of geometric hyperplanes of $S_k(3)$ under the action of the automorphism group $(\sigma_4\times \sigma_4\times \dots\times \sigma_4)\rtimes \sigma_k$.  To take care of the latter fact, it suffices to analyse behavior of Veldkamp lines of $S_{k-1}(3)$ under the action of $(\sigma_4\times \sigma_4\times \dots\times \sigma_4)\rtimes \sigma_{k-1}$. To this end in view,  we first compute the orbits of Veldkamp lines of $S_{k-1}(3)$ under the action of the small group $\sigma_4\times \sigma_4\times \dots\times \sigma_4$ and then group together the orbits that have elements in common when we act by $\sigma_{k-1}$ on $S_{k-1}(3)$. That is what Algorithm \ref{lst:algo3} does. In the algorithm, the action of the group $\sigma_4\times \sigma_4\times \dots\times \sigma_4$ is called a coordinates' permutation, while the action of $\sigma_{k-1}$ is called a dimension permutation.

\lstinputalgo[%
	label={lst:algo3},%
	caption={Veldkamp lines separation by orbits}%
]{veldkamp_lines_separation_by_orbits.algo}

%\newpage
\section{Veldkamp spaces of small ternary Segre varieties}\label{sec4}
\subsection{Veldkamp space of $S_2(3)$}
As the $k=1$ case is rather trivial ($\mathcal{V}(S_1(3)) \sim {\rm PG}(1,3)$), we shall pass directly to the next one. 
The Segre variety $S_2(3)$, which can be visualized as a $4 \times 4$-grid (see Figure 1),  features only two different types of geometric hyperplanes, namely singular hyperplanes and ovoids. A singular geometric hyperplane (denoted as $H_1$)  consists of two intersecting lines and there are 16 of them. An ovoid (denoted as $H_2$) is a set of four pairwise non-collinear points; there are 24 of them, explicitly listed in Table 1 employing a particular numbering of points of the variety depicted in Figure 1. $\mathcal{V}(S_2(3))$ possesses five distinct types of ordinary lines, illustrated in Figure 2, whose basic properties are listed in Table 2; in addition, Table 3 shows a detailed composition of all ovoidal lines. The two different types of trivial Veldkamp lines of $S_2(3)$ are portrayed in Figure 3.
\begin{figure}[pth!]
	\centering
	%\vspace*{1.0cm}
	\includegraphics[width=2.5truecm]{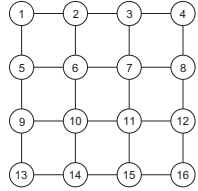}
	\caption{The adopted labeling of the points of $S_2(3)$.}
	\label{fig1}
\end{figure}
\begin{figure}[h]
	\centering
	%\vspace*{1.0cm}
	\centerline{\includegraphics[width=8.0truecm]{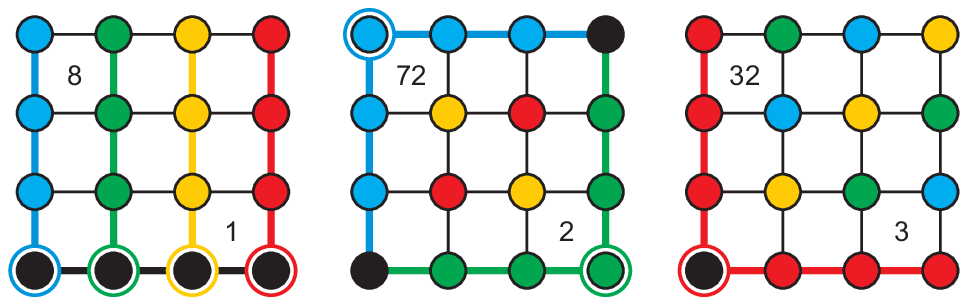}\hspace*{0.5cm}\includegraphics[width=5.0truecm]{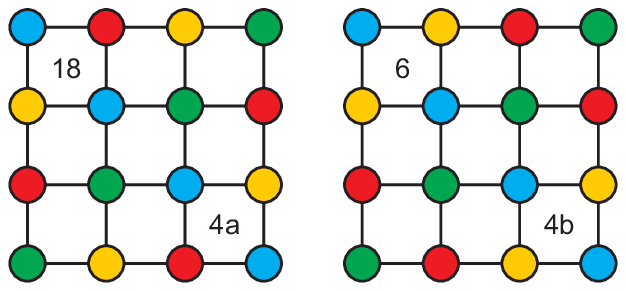}}
	\caption{A diagrammatic illustration of representatives of the five types of ordinary Veldkamp lines of $S_2(3)$. 
	In each representative the four hyperplanes are distinguished by different coloring, their core elements are depicted in black and deep points are encircled. The numbers in each subfigure stand for the type (bottom right) and its cardinality (top left).
	(Here types `4'and `4$^{\star}$' are denoted as `$4a$' and `$4b$', respectively.)}
	\label{fig2}
\end{figure}

\begin{figure}[t]
	\centering
	%\vspace*{1.0cm}
	\includegraphics[width=5.0truecm]{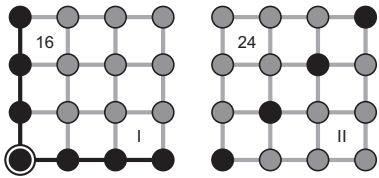} 
	\caption{A diagrammatic illustration of representatives of the two types of trivial Veldkamp lines of $S_2(3)$; a trivial line of type I/II comprises the $S_2(3)$ (gray) and a singular hyperplane/ovoid (black) counted three times. }
	\label{fig3}
\end{figure}

\begin{table}[pth!]
\begin{center}
\caption{The composition (rows) of 24 ovoids of $S_2(3)$. The `+' symbol indicates which point of $S_2(3)$ -- in  terms of the labeling of the latter given in Figure 1 -- lies in a given ovoid; for example, ovoid number 1 consists of points 4, 7, 10 and 13. The ovoids are grouped into six pairwise disjoint sets of four members, each such set representing a non-projective Veldkamp line of $S_2(3)$.} 
\vspace*{0.4cm} 
%\resizebox{\columnwidth}{!}{%
\begin{tabular}{||r|cccccccccccccccc||}
\hline \hline
     &  1  &  2  &  3  &  4  &  5  &  6  &  7  &  8  &  9  &  10    &   11   &   12   &  13    &  14    &  15    &  16  \\
\hline 
~1~  &     &     &     &  +  &     &     &  +  &     &     &  +     &        &        &   +    &        &        &          \\
~2~  &  +  &     &     &     &     &  +  &     &     &     &        &   +    &        &        &        &        &   +      \\
~3~  &     &     &  +  &     &     &     &     &  +  &  +  &        &        &        &        &   +    &        &          \\
~4~  &     &  +  &     &     &  +  &     &     &     &     &        &        &   +    &        &        &   +    &         \\
\hline
~5~  &     &     &  +  &     &     &     &     &  +  &     &  +     &        &        &   +    &        &        &         \\
~6~  &     &  +  &     &     &  +  &     &     &     &     &        &   +    &        &        &        &        &    +     \\
~7~  &     &     &     &  +  &     &     &  +  &     &  +  &        &        &        &        &   +    &        &         \\
~8~  &  +  &     &     &     &     &  +  &     &     &     &        &        &   +    &        &        &    +   &         \\
\hline
~9~  &     &     &     &  +  &     &  +  &     &     &     &        &   +    &        &   +    &        &        &         \\
~10~ &  +  &     &     &     &     &     &  +  &     &     &  +     &        &        &        &        &        &    +     \\
~11~ &     &     &  +  &     &  +  &     &     &     &     &        &        &   +    &        &    +   &        &          \\
~12~ &     &  +  &     &     &     &     &     &  +  &  +  &        &        &        &        &        &   +    &         \\
\hline
~13~ &     &  +  &     &     &     &     &  +  &     &     &        &        &   +    &   +    &        &        &          \\
~14~ &     &     &  +  &     &     &  +  &     &     &  +  &        &        &        &        &        &        &    +     \\
~15~ &     &     &     &  +  &  +  &     &     &     &     &  +     &        &        &        &        &   +    &          \\
~16~ &  +  &     &     &     &     &     &     &  +  &     &        &   +    &        &        &    +   &        &          \\
\hline
~17~ &     &  +  &     &     &     &     &     &  +  &     &        &   +    &        &   +    &        &        &          \\
~18~ &     &     &  +  &     &  +  &     &     &     &     &  +     &        &        &        &        &        &    +     \\
~19~ &     &     &     &  +  &     &  +  &     &     &  +  &        &        &        &        &        &    +   &          \\
~20~ &  +  &     &     &     &     &     &  +  &     &     &        &        &   +    &        &    +   &        &          \\
\hline
~21~ &     &     &  +  &     &     &  +  &     &     &     &        &        &   +    &    +   &        &        &          \\
~22~ &     &  +  &     &     &     &     &  +  &     &  +  &        &        &        &        &        &        &    +     \\
~23~ &     &     &     &  +  &  +  &     &     &     &     &        &   +    &        &        &    +   &        &         \\
~24~ &  +  &     &     &     &     &     &     &  +  &     &  +     &        &        &        &        &    +   &          \\
 \hline \hline
\end{tabular}
\end{center}
\label{tab1}
\end{table} 
%\clearpage
\begin{table}[pth!]
\begin{center}
\caption{The types of ordinary Veldkamp lines of $S_2(3)$. The first column gives the type, the next two columns tell us about how many points and lines belong to all the four geometric hyperplanes a line of the given type consists of, then we learn about the line's composition and, finally, the last column lists cardinalities for each type. The type denoted by a star is not projective.} 
\vspace*{0.2cm}
{\begin{tabular}{|l|c|c|c|c|c|} \hline \hline
\multicolumn{1}{|c|}{}  &  \multicolumn{2}{|c|}{} & \multicolumn{2}{|c|}{}  &  \multicolumn{1}{|c|}{}\\
%\cline{4-11}
\multicolumn{1}{|c|}{} &   \multicolumn{2}{|c|}{Core} & \multicolumn{2}{|c|}{Comp'n}  &\multicolumn{1}{|c|}{} \\
%\multicolumn{1}{||c|}{Hyperplane} & \multicolumn{1}{|c|}{Pts} & \multicolumn{1}{|c|}{Lns}  &  \multicolumn{1}{|c|}{0} & \multicolumn{1}{|c|}{1} & \multicolumn{2}{|c|}{0} & \multicolumn{1}{|c|}{3}
%& \multicolumn{1}{|c|}{deep} & \multicolumn{1}{|c|}{sing} & \multicolumn{1}{|c|}{ovoid} & \multicolumn{1}{|c|}{subq} &\multicolumn{1}{|c||}{} \\
 \cline{2-5}
Tp   & Ps   & Ls   & $H_1$ & $H_2$ & Crd  \\
\hline
1    &  4   & 1    & 4     & --    &  8 \\
\hline
2    &  2   & 0    & 2     & 2     &  72 \\
\hline
3    &  1   & 0    & 1     & 3     &  32  \\
\hline
4 &  0   & 0    & --    & 4     &  18 \\
\hline
$4^{\star}$ &  0   & 0    & --    & 4     &  6 \\
\hline \hline
\end{tabular}}
\end{center}
\label{tab2}
\end{table}

\begin{table}[pth!]
\begin{center}
\caption{The composition (rows) of 24 ovoidal Veldkamp lines of $S_2(3)$. The `+' symbol indicates which ovoid of $S_2(3)$ -- in terms of the labeling introduced in Table 1 -- belongs to a given Veldkamp line; for example, Veldkamp line number 1 consists of ovoids 1, 2, 11 and 12. The lines are grouped into four pairwise disjoint sets of six members, each such set representing a partition of the set of ovoids; the last six lines are not projective.} 
\vspace*{0.4cm} 
\resizebox{\columnwidth}{!}{%
\begin{tabular}{||r|cccccccccccccccccccccccc||}
\hline \hline
     & 1 & 2 & 3 & 4 & 5 & 6 & 7 & 8 & 9 & 10 & 11 & 12 & 13 & 14 & 15 & 16 & 17 & 18 & 19 & 20 & 21 & 22 & 23 & 24 \\
\hline 
~1~  & + & + &   &   &   &   &   &   &   &    &  + &  + &    &    &    &    &    &    &    &    &    &    &    &    \\
~2~  &   &   & + & + &   &   &   &   & + &  + &    &    &    &    &    &    &    &    &    &    &    &    &    &    \\
~3~  &   &   &   &   & + & + &   &   &   &    &    &    &    &    &    &    &    &    &  + & +  &    &    &    &    \\
~4~  &   &   &   &   &   &   & + & + &   &    &    &    &    &    &    &    &  + &  + &    &    &    &    &    &    \\
~5~  &   &   &   &   &   &   &   &   &   &    &    &    &  + &  + &    &    &    &    &    &    &    &    &  + & +  \\
~6~  &   &   &   &   &   &   &   &   &   &    &    &    &    &    &  + & +  &    &    &    &    &  + &  + &    &    \\
\hline
~7~  & + &   & + &   &   & + &   & + &   &    &    &    &    &    &    &    &    &    &    &    &    &    &    &    \\
~8~  &   & + &   & + & + &   & + &   &   &    &    &    &    &    &    &    &    &    &    &    &    &    &    &    \\
~9~  &   &   &   &   &   &   &   &   & + &    &  + &    &    &    &    &    &    &    &    &    &    &  + &    & +  \\
~10~ &   &   &   &   &   &   &   &   &   & +  &    &  + &    &    &    &    &    &    &    &    &  + &    &  + &    \\
~11~ &   &   &   &   &   &   &   &   &   &    &    &    &  + &    &    & +  &    &  + &  + &    &    &    &    &    \\
~12~ &   &   &   &   &   &   &   &   &   &    &    &    &    &  + &  + &    &  + &    &    &  + &    &    &    &    \\
\hline
~13~ & + &   &   & + &   &   &   &   &   &    &    &    &    &  + &    &  + &    &    &    &    &    &    &    &    \\
~14~ &   & + & + &   &   &   &   &   &   &    &    &    &  + &    &  + &    &    &    &    &    &    &    &    &    \\
~15~ &   &   &   &   & + &   &   & + &   &    &    &    &    &    &    &    &    &    &    &    &    &  + &  + &    \\
~16~ &   &   &   &   &   & + & + &   &   &    &    &    &    &    &    &    &    &    &    &    &  + &    &    & +  \\
~17~ &   &   &   &   &   &   &   &   & + &    &    &  + &    &    &    &    &    &  + &    & +  &    &    &    &    \\
~18~ &   &   &   &   &   &   &   &   &   & +  &  + &    &    &    &    &    &  + &    &  + &    &    &    &    &    \\
\hline
~19~ & + & + & + & + &   &   &   &   &   &    &    &    &    &    &    &    &    &    &    &    &    &    &    &    \\
~20~ &   &   &   &   & + & + & + & + &   &    &    &    &    &    &    &    &    &    &    &    &    &    &    &    \\
~21~ &   &   &   &   &   &   &   &   & + &  + &  + &  + &    &    &    &    &    &    &    &    &    &    &    &    \\
~22~ &   &   &   &   &   &   &   &   &   &    &    &    &  + & +  &  + &  + &    &    &    &    &    &    &    &    \\
~23~ &   &   &   &   &   &   &   &   &   &    &    &    &    &    &    &    &  + &  + &  + &  + &    &    &    &    \\
~24~ &   &   &   &   &   &   &   &   &   &    &    &    &    &    &    &    &    &    &    &    &  + &  + & +  & +  \\
 \hline \hline
\end{tabular}%
}
\end{center}
\label{tab3}
\end{table} 
%\clearpage

Now, the projective space PG$(3,3)$ features $(3^4 - 1)/(3-1) = 40$ points/planes and ${40 \choose 2}/{4 \choose 2} = 130$ lines.
One can verify that all 40 geometric hyperplanes and the total of 130 Veldkamp lines of types 1, 2, 3 and $4$ of $S_2(3)$ indeed form a space isomorphic to PG$(3,3)$\footnote{In the dual picture, one views $S_2(3)$ embedded as a hyperbolic quadric ${\cal Q}^+(3,3)$ in PG$(3,3)$ and checks that any of its 40 geometric hyperplanes originates from the intersection of this ${\cal Q}^+(3,3)$ with a certain plane of the PG$(3,3)$. }; the remaining six Veldkamp lines of type $4^{\star}$ have thus no counterparts in this PG$(3,3)$ and in what follows they will be referred to as {\it non-projective}. Hence, PG$(3,3) \subset \mathcal{V}(S_2(3))$, which contrasts PG$(3,2) \sim \mathcal{V}(S_2(2))$ \cite{seg2}. 
%$S_2(3)$ features two distinct types of trivial Veldkamp lines, whose representatives are depicted in Figure 3.

\subsection{Veldkamp space of $S_3(3)$}

Next, for each type of the line of $\mathcal{V}(S_2(3))$, as well as for either type of the two trivial, let us take one representative and blow it up in order to readily see that $\mathcal{V}(S_3(3))$ exhibits six different types of points, whose principal properties are described in Table 4 and whose representatives are depicted in Figure 4. The relevant projective space is now PG$(7,3)$, which is endowed with $ (3^8 -1)/2 = 3280$ points/hyperplanes and ${3280 \choose 2}/{4 \choose 2} = 896\,260$ lines. 
From Table 4 we infer that there are 3280 geometric hyperplanes of $S_3(3)$ of types $H_1$, $H_2$, $H_3$, $H_4$ and $H_5$, which will be the points of PG$(7,3) \subset \mathcal{V}(S_3(3))$; the remaining 144 geometric hyperplanes of type $H_5^{\star}$ have no counterparts in this PG$(7,3)$, since each of them is the blow-up of a non-projective Veldkamp line of $S_2(3)$.

\begin{table}[t]
\begin{center}
\caption{The six types of geometric hyperplanes of the Segre variety $S_3(3)$. The first column gives the type (`Tp') of a hyperplane, which is followed by the number of points (`Pts') and lines (`Lns') it contains, and the number of points of given order. The next three columns tell us about how many of 12 $S_2(3)$'s are fully located  in the hyperplane ($D$) and/or share with it a hyperplane of type $H_1$ (i.\,e. a singular hyperplane) or $H_2$ (i.\,e. an ovoid). The VL-column lists the types of (both ordinary and trivial) Veldkamp lines of $S_2(3)$  we get by projecting a hyperplane of the given type into an $S_2(3)$ along the lines of all three distinguished spreads. Finally, for each hyperplane type we give its cardinality (`Crd'), the corresponding large orbit  of $2 \times 2 \times 2$ arrays over $GF(3)$ (`BS') taken from Table 4 of \cite{brst}, and its weight, or rank in the language of \cite{brst} (`W'). The type denoted by a star symbol is not projective.} 
\vspace*{0.3cm}
{\begin{tabular}{|l|c|c|c|c|c|c|c|c|c|l|r|c|c|} \hline \hline
%\multicolumn{1}{|c|}{} & \multicolumn{1}{|c|}{} & \multicolumn{1}{|c|}{}  &  \multicolumn{4}{|c|}{} & \multicolumn{3}{|c|}{}  &\multicolumn{1}{|c|}{} &\multicolumn{1}{|c|}{} &\multicolumn{1}{|c|}{}  & \multicolumn{1}{|c|}{}\\
%\cline{4-11}
\multicolumn{1}{|c|}{} & \multicolumn{1}{|c|}{} &
\multicolumn{1}{|c|}{}  &  \multicolumn{4}{|c|}{Points of
Order} & \multicolumn{3}{|c|}{$S_2(3)$'s of Type} & \multicolumn{1}{|c|}{} & \multicolumn{1}{|c|}{} &\multicolumn{1}{|c|}{} 
& \multicolumn{1}{|c|}{}\\
 \cline{4-10}
Tp & Pts & Lns  & ~0~ & ~1~ & 2 & ~3~ & ~D~ & $H_1$ & $H_2$ & VL & Crd & BS & W \\
\hline
1 & 37 & 21  & 0 & 0 & 27   & 10 & 3 & 9 & 0 & I      & 64  & 2  & 1 \\
\hline
2 & 28 & 12  & 0 & 12 & 12  & 4 & 1 & 8 & 3 & II,\,1 & 288  & 3  & 2 \\
\hline
3 & 22 & 6   & 4 & 12 & 6  & 0 & 0 & 6 & 6 &  2     & 1728 & 4  & 2 \\
\hline
4 & 19 & 3   & 9 & 9 & 0    & 1 & 0 & 3 & 9 &  3     & 768  & 5  & 3 \\
\hline
5 & 16  & 0   & 16 & 0 & 0    & 0 & 0 & 0 & 12 &  4     & 432  & 6  & 3 \\
\hline
$5^{\star}$ & 16  & 0   & 16 & 0 & 0    & 0 & 0 & 0 & 12 &  4$^{\star}$    & 144  & $-$  & $-$ \\
\hline \hline
\end{tabular}}
\end{center}
\label{tab4}
\end{table}

\begin{figure}[h]
	\centering
	%\vspace*{1.0cm}
	\includegraphics[width=6.5truecm]{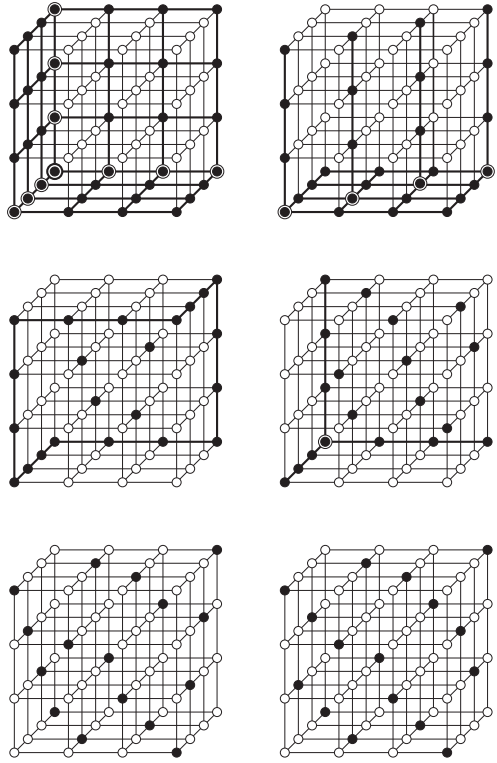} 
	\caption{A diagrammatic illustration of representatives of the six distinct types (numbered consecutively from top left ($H_1$, a singular hyperplane) to bottom right ($H_5^{\star}$, a non-projective ovoid) of geometric hyperplanes of $S_3(3)$, where the latter is represented by a $4 \times 4 \times 4$-grid. As in Figures 2 and 3, deep points are encircled.}
	\label{fig4}
\end{figure}

%\subsubsection{Two kinds of ovoids of $S_3(3)$: A closer look}
The existence of two different types of ovoids (projective and non-projective) of $S_3(3)$ stems, obviously, from the existence of two different kinds of ovoidal Veldkamp lines of $S_2(3)$.  But there is a more intricate distinction between the two types when some graph theory is invoked -- namely the celebrated {\it Dyck} graph \cite{dyck}. This graph is of valency three and has 32 vertices, 48 edges, 12 octagons, being of girth six, diameter five and having chromatic number equal to two; it is a non-planar graph of minimal genus $g=1$ (i.e., it can be embedded on the torus). Now, let us take a disjoint union of two projective ovoids and consider the graph whose vertices are the points of the ovoids, with two vertices being connected by an edge if they
are collinear; the graph we get is isomorphic to {\it either} the Dyck graph {\it or} a disjoint union of four cubical graphs. If, however, we take a disjoint union of two non-projective ovoids, the graph we obtain is {\it never} isomorphic to the Dyck graph. 
Moreover, taking the symmetric difference of two non-projective ovoids that share four points situated at maximum distance from each other, we get the so-called {\it Nauru} graph \cite{epst} -- a symmetric bipartite cubic graph with 24 vertices and 36 edges isomorphic to  the generalized Petersen graph $G(12,5)$.

%\subsection{Veldkamp lines of $S_3(3)$}
Disregarding first all geometric hyperplanes of type $H_5^{\star}$ and using only geometric and combinatorial arguments, we have found altogether 66 different types of Veldkamp lines of $S_3(3)$. Four of them, having 2268 elements in total, are not projective and the remaining 62 types, amounting to 896\,260 members, were all found to belong to the PG$(7,3)$. Table 5 gives a summary of basic combinatorial characteristics for each type.
As one sees, there are a number of cases where two (and in one case even three) types are characterized by the same string of parameters. In these cases it was necessary to look in more detail at the properties of the core to see the difference. In particular,
\begin{itemize}
\item 16 vs 17: the two lines of the core are, respectively, concurrent or skew.
\item 21 vs 22: out of four points of order zero of both $H_3$'s, one or none belongs to the core, respectively.
\item 27 vs 28 vs 29: out of the four points of order zero of each $H_3$, in type 27 one belongs to the core, whereas in type 28 two enjoy this property; concerning type 29,  for two $H_3$'s there are two zero-order points that lie in the core and for either of the other two $H_3$'s it is only one such point in the core. 
\item 31 vs 32: out of the four points of order zero of both $H_3$'s, two or three belong to the core, respectively.
\item 33 vs 34: in the former case all the four zero-order points of $H_3$ belong to the core, whilst in the latter case it is only two of them; also, in type 33 the two deep points of $H_4$'s  are at distance two, whilst in type 34 they are at distance three. 
\item 39 vs 40: the single deep point of $H_4$ does or does not belong to the core, respectively.
\item 44 vs 45: out of the four points of order zero of each $H_3$, none or several belong to the core, respectively.
\item 46 vs 47: the deep points of $H_4$'s are at distance two or three, respectively.
\item 48 vs 49: the deep point of the $H_4$ does not or does belong to the core, respectively.
\item 53 vs 54: out of the four points of the core, two coincide with two points of order two in either of the two $H_3$'s, these two points being at distance two or three, respectively.
\item 55 vs 56: the deep points of $H_4$'s are at distance two or three, respectively; also, in the latter case one of the two deep points belongs to the core.
\item 57 vs 58: none or all the four deep points of $H_4$'s belong to the core, respectively.
\end{itemize}

\begin{table}[h]
\begin{center}
\caption{The 66 types of ordinary Veldkamp lines of $S_3(3)$ that can be discerned from each other on solely geometrical grounds, i.\,e. without invoking group theoretical arguments; the four starred ones are not projective. The notation follows that of Table 2. (Compare with Table 4 of \cite{seg2} or Table 2 of \cite{gs}.)} \vspace*{0.5cm}
{\begin{tabular}{|l|c|c|c|c|c|c|c|r|} \hline \hline
%\multicolumn{1}{|c|}{}  &  \multicolumn{2}{|c|}{} & \multicolumn{5}{|c|}{}  &  \multicolumn{1}{|c|}{}\\
%\cline{4-11}
\multicolumn{1}{|c|}{} &   \multicolumn{2}{|c|}{Core} & \multicolumn{5}{|c|}{Composition}  &\multicolumn{1}{|c|}{} \\
%\multicolumn{1}{||c|}{Hyperplane} & \multicolumn{1}{|c|}{Pts} & \multicolumn{1}{|c|}{Lns}  &  \multicolumn{1}{|c|}{0} & \multicolumn{1}{|c|}{1} & \multicolumn{2}{|c|}{0} & \multicolumn{1}{|c|}{3}
%& \multicolumn{1}{|c|}{deep} & \multicolumn{1}{|c|}{sing} & \multicolumn{1}{|c|}{ovoid} & \multicolumn{1}{|c|}{subq} &\multicolumn{1}{|c||}{} \\
 \cline{2-8}
Tp & ~Pts & ~Lns  & $H_1$ & $H_2$ & $H_3$ & $H_4$ & $H_5$ & Crd  \\
\hline
~1  & ~28 & ~15  & 4  & -- & -- & -- & -- &  48 \\
\hline
~2  & ~22 & ~10  & 2  & 2  & -- & -- & -- &  864 \\
\hline
~3  & ~19 & ~9   & 1  & 3  & -- & -- & -- & 384  \\
\hline
~4  & ~18 & ~6   & 2  & -- & 2  & -- & -- &  864 \\
\hline
~5  & ~16 & ~8   & -- &  4 & -- & -- & -- & 216  \\
~5$^{\star}$  & ~16 & ~8   & -- &  4 & -- & -- & -- & 72  \\
\hline
~6  & ~16 & ~4   & -- & 4  & -- & -- & -- & 72  \\
\hline
~7  & ~15 & ~4   & 1  & 1  & 2  & -- & -- &  10368 \\
\hline
~8  & ~13 & ~4   & 1  & -- & 3  & -- & -- &  3456 \\
\hline
~9  & ~13 & ~3   & 1  & 1  & -- & 2  & -- &  3456 \\
10  & ~13 & ~3   & 1  & -- & 3  & -- & -- &  6912 \\
11  & ~13 & ~3   & -- & 3  & -- & 1  & -- &  2304 \\
\hline
12  & ~12 & ~3   & -- & 2  & 2  & -- & -- & 20736 \\
\hline
13  & ~12 & ~2   & 1  & -- & 2  & 1  & -- & 20736 \\
14  & ~12 & ~2   & -- & 2  & 2  & -- & -- &  2592 \\
\hline
15  & ~10 & ~3   & 1  & -- & -- & 3  & -- & 256  \\
\hline
16  & ~10 & ~2   & -- & 1  & 3  & -- & -- & 20736  \\
17  & ~10 & ~2   & -- & 1  & 3  & -- & -- & 3456  \\
\hline
18  & ~10 & ~1   & -- & 2  & -- & 2  & -- & 3456  \\
19  & ~10 & ~1   & -- & 1  & 3  & -- & -- & 13824 \\
\hline
20  & ~10 & ~0   & 1  & -- & 1  & 1  & 1  & 6912  \\
\hline
21  &  ~9 & ~2   & -- &  1 & 2  & 1  & -- & 20736  \\
22  &  ~9 & ~2   & -- &  1 & 2  & 1  & -- & 20736  \\
\hline
23  &  ~9 & ~1   & -- & 1  & 2  & 1  & -- & 20736 \\
\hline
24  &  ~9 & ~0   & 1  & -- & 1  & -- & 2  & 6912  \\
25  &  ~9 & ~0   & 1  & -- & -- & 2  & 1  & 6912  \\
\hline
26  &  ~8 & ~2   & -- & -- & 4  & -- & -- & 3888 \\
26$^{\star}$  &  ~8 & ~2   & -- & -- & 4  & -- & -- & 1296 \\
\hline
27&  ~8 & ~1   & -- & -- & 4  & -- & -- & 10368  \\
28&  ~8 & ~1   & -- & -- & 4  & -- & -- & 10368  \\
29&  ~8 & ~1   & -- & -- & 4  & -- & -- & 41472 \\
\hline
30  &  ~8 & ~0   & -- & 2  & -- & -- & 2  & 3888  \\
31&  ~8 & ~0   & -- & 1  & 2  & -- & 1  & 10368  \\
32&  ~8 & ~0   & -- & 1  & 2  & -- & 1  & 20736  \\
33&  ~8 & ~0   & -- & 1  & 1  & 2  & -- & 10368  \\
34&  ~8 & ~0   & -- & 1  & 1  & 2  & -- & 20736  \\
35  &  ~8 & ~0   & -- & -- & 4  & -- & -- & 5184 \\
\hline
36  &  ~7 & ~1   & -- & 1  & -- & 3  & -- & 2304  \\
37  &  ~7 & ~1   & -- & -- & 3  & 1  & -- & 41472 \\
\hline
38  &  ~7 & ~0   & -- & 1  & 1  & 1  & 1  & 41472 \\
39  &  ~7 & ~0   & -- & -- & 3  & 1  & -- & 2304  \\
40  &  ~7 & ~0   & -- & -- & 3  & 1  & -- & 69120  \\
 \hline \hline
\end{tabular}}
\end{center}
\label{tab5}
\end{table}
\clearpage
\addtocounter{table}{-1}
\vspace*{-0.3cm}
\begin{table}[h]
\begin{center}
\caption{Continued.} \vspace*{0.5cm}
{\begin{tabular}{|l|c|c|c|c|c|c|c|r|} \hline \hline
%\multicolumn{1}{|c|}{}  &  \multicolumn{2}{|c|}{} & \multicolumn{5}{|c|}{}  &  \multicolumn{1}{|c|}{}\\
%\cline{4-11}
\multicolumn{1}{|c|}{} &   \multicolumn{2}{|c|}{Core} & \multicolumn{5}{|c|}{Composition}  &\multicolumn{1}{|c|}{} \\
%\multicolumn{1}{||c|}{Hyperplane} & \multicolumn{1}{|c|}{Pts} & \multicolumn{1}{|c|}{Lns}  &  \multicolumn{1}{|c|}{0} & \multicolumn{1}{|c|}{1} & \multicolumn{2}{|c|}{0} & \multicolumn{1}{|c|}{3}
%& \multicolumn{1}{|c|}{deep} & \multicolumn{1}{|c|}{sing} & \multicolumn{1}{|c|}{ovoid} & \multicolumn{1}{|c|}{subq} &\multicolumn{1}{|c||}{} \\
 \cline{2-8}
Tp & ~Pts & ~Lns  & $H_1$ & $H_2$ & $H_3$ & $H_4$ & $H_5$ & Crd  \\
\hline
41  &  ~6 & ~1   & -- & -- & 2  & 2  & -- & 20736  \\
\hline
42  &  ~6 & ~0   & -- & 1  & 1  & -- & 2  & 10368 \\
43  &  ~6 & ~0   & -- & 1  & -- & 2  & 1  & 10368  \\
44&  ~6 & ~0   & -- & -- & 3  & -- & 1  & 6912  \\
45&  ~6 & ~0   & -- & -- & 3  & -- & 1  & 62208 \\
46&  ~6 & ~0   & -- & -- & 2  & 2  & -- & 62208 \\
47&  ~6 & ~0   & -- & -- & 2  & 2  & -- & 6912 \\
\hline
48&  ~5 & ~0   & -- & -- & 2  & 1  & 1  & 82944  \\
49&  ~5 & ~0   & -- & -- & 2  & 1  & 1  & 20736  \\
50  &  ~5 & ~0   & -- & -- & 1  & 3  & --  & 20736 \\
\hline
51  &  ~4 & ~1   & -- & -- & -- & 4  & -- & 1728  \\
51$^{\star}$  &  ~4 & ~1   & -- & -- & -- & 4  & -- & 576  \\
\hline
52  &  ~4 & ~0   & -- & 1  & -- & -- & 3  & 1728  \\
53&  ~4 & ~0   & -- & -- & 2  & -- & 2  & 20736  \\
54&  ~4 & ~0   & -- & -- & 2  & -- & 2  & 10368  \\
55&  ~4 & ~0   & -- & -- & 1  & 2  & 1  & 20736  \\
56&  ~4 & ~0   & -- & -- & 1  & 2  & 1  & 13824  \\
57&  ~4 & ~0   & -- & -- & -- & 4  & -- & 3456  \\
58&  ~4 & ~0   & -- & -- & -- & 4  & -- & 576  \\
\hline
59   &  ~3 & ~0   & -- & -- & 1  & 1  & 2  & 13824  \\
\hline
60   &  ~2 & ~0   & -- & -- & -- & 2  & 2  & 10368 \\
\hline
61   &  ~1 & ~0   & -- & -- & -- & 1  & 3  & 2304  \\
\hline
62   &  ~0 & ~0   & -- & -- & -- & -- & 4  & 756  \\
62$^{\star}$   &  ~0 & ~0   & -- & -- & -- & -- & 4  & 324  \\
 \hline \hline
\end{tabular}}
\end{center}
\label{tab5}
\end{table} 

Table 5 reveals a number of interesting properties of the Veldkamp lines of $S_3(3)$. First, we notice that there are only two types (62 and 62$^{\star}$) whose core is an empty set; all the four geometric hyperplanes of this line are of the same kind, namely ovoids. Next, there are 16 different types of Veldkamp lines each of which consists of  hyperplanes of the same type. Explicitly, these are types 1 ($H_1$), 5, 
5$^{\star}$, 6 ($H_2$), 26, 26$^{\star}$, 27, 28, 29, 35 ($H_3$), 51, 51$^{\star}$, 57, 58 ($H_4$) and 62 and 62$^{\star}$ ($H_5$); note that they all have an even number of points in their cores. On the
other hand, we find only two  types (compared with as many as nine in the binary case) where all constituting hyperplanes are of different types; these are Veldkamp lines of types 20 and 38. The most ``abundant" 
type of geometric hyperplanes is $H_3$, which occurs in 43 types of lines, followed by $H_4$ (33 types), $H_2$ (27 types), $H_5$ (23 types) and, finally, by $H_1$ (13 types). Remarkably, $H_3$ also prevails in 
multiplicity $\geq 3$ (16 types), whereas in ``singles'' the primacy belongs to $H_2$ (17 types). Further, one observes that if a line contains two $H_1$'s or two $H_2$'s, its remaining two hyperplanes are also of 
the same type. One also notes that the number of types of Veldkamp lines that feature just three $H_i$'s, $i = 1, 2,\ldots, 5$, is zero, two, ten, three and two, respectively, and that there is no non-projective Veldkamp line endowed with $H_1$'s.

The last two properties clearly indicate that geometric hyperplanes of type one play, as in the binary case, a distinguished role. In particular, they enable us to define the {\it weight} (or, in the language of Bremner and Stavrou \cite{brst}, the {\it rank}) of a geometric hyperplane in an inductive way as follows. Let us take any geometric hyperplane of type one to be of weight one. We call a geometric hyperplane to be of weight two if it is found to lie on a line of $\mathcal{V}(S_3(3))$ defined by two distinct $H_1$'s. From Table 5 we infer 
that apart from type one (consisting solely of $H_1$'s), there are only two more types of such a line; type two, where the two remaining hyperplanes are
$H_2$'s and type four, featuring two $H_3$'s. Hence, the geometric hyperplanes of weight two are those of type two and three. Next, we
call a geometric hyperplane to be of weight three if it is found to lie on a line of $\mathcal{V}(S_3(3))$ defined by an $H_1$ and any hyperplane of weight two; from Table 5 we see that both $H_4$'s and $H_5$'s are of weight three (see also Table 4). Thus, weight three is the highest one for geometric hyperplanes of $S_3(3)$.

Concerning non-projective lines of $\mathcal{V}(S_3(3))$, we have already seen four types of them consisting solely of projective points of $\mathcal{V}(S_3(3))$. Any other non-projective line necessarily entails a non-projective ovoid ($H_5^{\star}$). There are a variety of them, many featuring just two geometric hyperplanes. Those having size four fall into five distinct types, as shown in Table 6; note that $H_2$'s and $H_3$'s are the only non-ovoidal types involved. As per lines of size two, it is worth mentioning a line defined by two non-projective ovoids sharing four points that are at maximum mutual distance from each other, or a line defined by a non-projective ovoid and an $H_4$ having a single point in common, this being the deep point of the latter. 

\begin{table}[h]
\begin{center}
\caption{The five types of non-projective Veldkamp lines of $S_3(3)$ that contain non-projective ovoids.} \vspace*{0.5cm}
{\begin{tabular}{|l|c|c|c|c|c|c|c|c|r|} \hline \hline
%\multicolumn{1}{|c|}{}  &  \multicolumn{2}{|c|}{} & \multicolumn{5}{|c|}{}  &  \multicolumn{1}{|c|}{}\\
%\cline{4-11}
\multicolumn{1}{|c|}{} &   \multicolumn{2}{|c|}{Core} & \multicolumn{6}{|c|}{Composition}  &\multicolumn{1}{|c|}{} \\
%\multicolumn{1}{||c|}{Hyperplane} & \multicolumn{1}{|c|}{Pts} & \multicolumn{1}{|c|}{Lns}  &  \multicolumn{1}{|c|}{0} & \multicolumn{1}{|c|}{1} & \multicolumn{2}{|c|}{0} & \multicolumn{1}{|c|}{3}
%& \multicolumn{1}{|c|}{deep} & \multicolumn{1}{|c|}{sing} & \multicolumn{1}{|c|}{ovoid} & \multicolumn{1}{|c|}{subq} &\multicolumn{1}{|c||}{} \\
 \cline{2-9}
Tp             & ~Pts & ~Lns  & $H_1$ & $H_2$ & $H_3$ & $H_4$ & $H_5$ & $H_5^{\star}$ & Crd  \\
\hline
30$^{\star}$   &  ~8  & ~0    & --    & 2     & --    & --    & --    &   2           & 1296  \\
\hline
44$^{\star}$   &  ~6  & ~0    & --    & --    & 3     & --    & --    &   1           & 2304 \\
\hline
52$^{\star}$   &  ~4  & ~0    & --    & 1     & --    & --    & --    &   3           &  576 \\
\hline
$62_1^{\star}$ &  ~0  & ~0    & --    & --    & --    & --    & 2     &   2           &  864 \\
$62_2^{\star}$ &  ~0  & ~0    & --    & --    & --    & --    & --    &   4           &  360  \\
 \hline \hline
\end{tabular}}
\end{center}
\label{tab6}
\end{table}

%\newpage
To conclude this subsection, one has to stress that if we take into account a fine structure of the automorphism group of $S_3(3)$ we can further refine our classification of projective Veldkamp lines of $S_3(3)$.
In particular, we find seven types in Table 5 such that each splits into two subtypes with cardinalities as shown in Table 7 below; we shall see the importance of these splittings when dealing with classification of geometric hyperplanes of $S_4(3)$. 

\bigskip
\begin{table}[h]
\begin{center}
\caption{The seven particular types of projective lines of $\mathcal{V}(S_3(3))$ that split into two subtypes each when analyzed in group-theoretical terms.} \vspace*{0.5cm}
{\begin{tabular}{|l|c|c|c|c|c|c|c|r|} \hline \hline
%\multicolumn{1}{|c|}{}  &  \multicolumn{2}{|c|}{} & \multicolumn{5}{|c|}{}  &  \multicolumn{1}{|c|}{}\\
%\cline{4-11}
\multicolumn{1}{|c|}{} &   \multicolumn{2}{|c|}{Core} & \multicolumn{5}{|c|}{Composition}  &\multicolumn{1}{|c|}{} \\
%\multicolumn{1}{||c|}{Hyperplane} & \multicolumn{1}{|c|}{Pts} & \multicolumn{1}{|c|}{Lns}  &  \multicolumn{1}{|c|}{0} & \multicolumn{1}{|c|}{1} & \multicolumn{2}{|c|}{0} & \multicolumn{1}{|c|}{3}
%& \multicolumn{1}{|c|}{deep} & \multicolumn{1}{|c|}{sing} & \multicolumn{1}{|c|}{ovoid} & \multicolumn{1}{|c|}{subq} &\multicolumn{1}{|c||}{} \\
 \cline{2-8}
Tp & ~Pts & ~Lns  & $H_1$ & $H_2$ & $H_3$ & $H_4$ & $H_5$ & Crd  \\
\hline 
26a  &  ~8 & ~2   & -- & -- & 4  & -- & -- & 1296 \\
26b  &  ~8 & ~2   & -- & -- & 4  & -- & -- & 2592 \\
\hline
30a  &  ~8 & ~0   & -- & 2  & -- & -- & 2  & 1296  \\
30b  &  ~8 & ~0   & -- & 2  & -- & -- & 2  & 2592  \\
\hline
40a  &  ~7 & ~0   & -- & -- & 3  & 1  & -- & 27648  \\
40b  &  ~7 & ~0   & -- & -- & 3  & 1  & -- & 41472  \\
\hline
45a  &  ~6 & ~0   & -- & -- & 3  & -- & 1  & 20736 \\
45b  &  ~6 & ~0   & -- & -- & 3  & -- & 1  & 41472 \\
46a  &  ~6 & ~0   & -- & -- & 2  & 2  & -- & 20736 \\
46b  &  ~6 & ~0   & -- & -- & 2  & 2  & -- & 41472 \\
\hline
48a  &  ~5 & ~0   & -- & -- & 2  & 1  & 1  & 41472  \\
48b  &  ~5 & ~0   & -- & -- & 2  & 1  & 1  & 41472  \\
\hline
62a  &  ~0 & ~0   & -- & -- & -- & -- & 4  & 108  \\
62b  &  ~0 & ~0   & -- & -- & -- & -- & 4  & 648  \\
 \hline \hline
\end{tabular}}
\end{center}
\label{tab7}
\end{table}

\begin{table}[pth!]
\begin{center}
\caption{Those 15 types of lines of $\mathcal{V}(S_3(3))$ that are extensions of lines of $\mathcal{V}(S_3(2))$, the types of the latter being explicitly indicated in the last column (adopting the numbering from Table 4 of \cite{seg2} or Table 2 of \cite{gs}).} \vspace*{0.5cm}
{\begin{tabular}{|l|c|c|c|c|c|c|c|c|c|} \hline \hline
%\multicolumn{1}{|c|}{}  &  \multicolumn{2}{|c|}{} & \multicolumn{5}{|c|}{}  &  \multicolumn{1}{|c|}{}\\
%\cline{4-11}
\multicolumn{1}{|c|}{} &   \multicolumn{2}{|c|}{Core} & \multicolumn{6}{|c|}{Composition}  &\multicolumn{1}{|c|}{} \\
%\multicolumn{1}{||c|}{Hyperplane} & \multicolumn{1}{|c|}{Pts} & \multicolumn{1}{|c|}{Lns}  &  \multicolumn{1}{|c|}{0} & \multicolumn{1}{|c|}{1} & \multicolumn{2}{|c|}{0} & \multicolumn{1}{|c|}{3}
%& \multicolumn{1}{|c|}{deep} & \multicolumn{1}{|c|}{sing} & \multicolumn{1}{|c|}{ovoid} & \multicolumn{1}{|c|}{subq} &\multicolumn{1}{|c||}{} \\
 \cline{2-9}
Tp             & ~Pts & ~Lns  & $H_1$ & $H_2$ & $H_3$ & $H_4$ & $H_5$ & $H_5^{\star}$ & BVL  \\
\hline
~1              &  28  & 15    &  4    & --    & --    & --    & --    &   --          &  1 \\
\hline
~2              &  22  & 10    &  2    & 2     & --    & --    & --    &   --          &  2 \\
\hline
~3              &  19  & 9     &  1    & 3     & --    & --    & --    &   --          &  6 \\
\hline
~4              &  18  & 6     &  2    & --    & 2     & --    & --    &   --          &  3 \\
\hline
~6              &  16  & 4     &  --   & 4     & --    & --    & --    &   --          &  10 \\
\hline
~7              &  15  & 4     &  1    & 1     & 2     & --    & --    &   --          &  5 \\
\hline
10             &  13  & 3     &  1    & --    & 3     & --    & --    &   --          &  9 \\
11             &  13  & 3     &  --   & 3     & --    & 1     & --    &   --          &  11 \\
\hline
15             &  10  & 3     &  1    & --    & --    & 3     & --    &   --          &  41 \\
\hline
16             &  10  & 2     &  --   & 1     & 3     & --    & --    &   --          &  16 \\
17             &  10  & 2     &  --   & 1     & 3     & --    & --    &   --          &  23 \\
\hline
18             &  10  & 1     &  --   & 2     & --    & 2     & --    &   --          &  26 \\
\hline
20             &  10  & 0     &  1    & --    & 1     & 1     & 1     &   --          &  20 \\
\hline
23             &  9   & 1     &  --   & 1     & 2     & 1     & --    &   --          &  19 \\
\hline
44$^{\star}$   &  6   & 0     &  --   & --    & 3     & --    & --    &   1           &  28 \\
 \hline \hline
\end{tabular}}
\end{center}
\label{tab8}
\end{table}

\subsection{Relating binary and ternary Veldkamp spaces}   
At this point we will make a slight digression from the main line of the paper in order to shed some unexpected light on the findings of the previous subsection.
The occurrence of both non-projective points and non-projective lines of $\mathcal{V}(S_3(3))$ is very interesting because its binary counterpart, $\mathcal{V}(S_3(2)) \sim $ PG$(7,2)$, is {\it fully} projective.
Now, $S_3(3)$ contains 64 different copies of $S_3(2)$. Let us pick up one such copy and choose in it a Veldkamp line. A natural question emerges: what type does this Veldkamp line have to be to admit extension into (i.\,e., to be a trace of) a Veldkamp line of the ambient $S_3(3)$? Or, rephrased differently, given a Veldkamp line of $S_3(3)$, is this still a Veldkamp line when restricted to one (or even several) of its $S_3(2)$'s and, if so, what type is it? 
To answer these questions, we proceeded as follows. First, we checked the extendability of geometric hyperplanes of $S_3(2)$ (first classified in \cite{seg2})  and found out that for $1 \leq s \leq 3$ $H_s$ of $S_3(2)$ is, respectively, extendible to $H_s$ of $S_3(3)$,  $H_5$ of $S_3(2)$ extends into $H_4$ of $S_3(3)$ and $H_4$ of $S_3(2)$ is not extendible. With this information in mind, we took from Figure 6 of  \cite{seg2} the diagrammatic portrayal of a representative of each type of Veldkamp line of $S_3(2)$ that does not contain $H_4$'s, `pasted' it into a copy of $S_3(2)$ selected in $S_3(3)$  drawn as an  $4 \times 4 \times 4 $-grid and tried to complete this Veldkamp line into a Veldkamp line of the ambient $S_3(3)$.  
The results of our analysis are given in Table 8. We see that, not surprisingly, it is only only 15 (out of 41) types of lines of  $\mathcal{V}(S_3(2))$ that can be extended to lines of $\mathcal{V}(S_3(3))$. Yet, what really comes as a big surprise is that one of them -- being, of course, projective --  extends into a line of $\mathcal{V}(S_3(3))$ that is {\it not} projective (see the last row of Table 8)! 
We are sure that this peculiar feature  is telling us something crucial about the fine structure of $\mathcal{V}(S_3(3))$ that at this stage cannot be fully appreciated.

%\newpage
\subsection{Geometric hyperplanes of $S_4(3)$}
Applying the blow-up recipe on a representative of each projective line type of $\mathcal{V}(S_3(3))$ given in Table 5 as well as of each type of trivial Veldkamp lines, we found $(3^{16} -1)/2$ = 
21\,523\,360 geometric hyperplanes of $S_4(3)$ that are in bijection with the points of PG$(15,3) \subset \mathcal{V}(S_4(3))$,  falling into 43 different types as listed in Table 9.

\begin{table}[pth!]
\centering
\caption{The 43 types of projective  geometric hyperplanes of $S_4(3)$ that can be discerned from each other on solely geometrical grounds; also shown is partition of hyperplane types into 20 classes according to the number of points/lines. As in Table 4, one first gives the type (`Tp') of a hyperplane, then the number of points (`Pts') and lines (`Lns') it contains, and the number of points of given order. The next six columns tell us about how many of 16 $S_3(3)$'s are fully located (`D') in the hyperplane and/or share with it a hyperplane of type $H_i$ (see Table 4). The VL-column lists the types of ordinary and/or trivial Veldkamp lines of $S_3(3)$  we get by projecting a hyperplane of the given type into $S_3(3)$'s along the lines of all four distinguished spreads. Finally, for each hyperplane type we give its cardinality (`Crd'), the corresponding large orbit  of $2 \times 2 \times 2 \times 2$ arrays over $GF(3)$ (`BS') taken from Table 6 of \cite{brst}, and its weight/rank (`W'). (Compare with Table 5 of \cite{seg2}.)}
\bigskip
\resizebox{\columnwidth}{!}{% 
{\begin{tabular}{|r|r|r|r|r|r|r|r|c|c|c|c|c|c|r|r|c|c|} \hline \hline
%\multicolumn{1}{|c|}{} & \multicolumn{1}{|c|}{} & \multicolumn{1}{|c|}{}  &  \multicolumn{5}{|c|}{}                        & \multicolumn{6}{|c|}{}                      & \multicolumn{1}{|c|}{}  & \multicolumn{1}{|c|}{} & \multicolumn{1}{|c|}{} & \multicolumn{1}{|c|}{} \\
%\cline{4-11}
\multicolumn{1}{|c|}{} & \multicolumn{1}{|c|}{} & \multicolumn{1}{|c|}{}  &  \multicolumn{5}{|c|}{$\#$ of Points of Order} & \multicolumn{6}{|c|}{$\#$ of $S_3(3)$'s of Type} & \multicolumn{1}{|c|}{}  & \multicolumn{1}{|c|}{} & \multicolumn{1}{|c|}{} & \multicolumn{1}{|c|}{} \\
 \cline{4-14}
Tp  &  Pts & Lns  &  0  &  1  &  2  &  3  &  4  &  D   & $H_1$  & $H_2$ & $H_3$ & $H_4$ &  $H_5$ & VL               & Crd       & BS    & W\\
\hline
~1  &  175 & 148  &  0  &  0  &  0  & 108 & 67  &  4   &   12   &   0   &   0   &   0   &   0    & I                & 256       & 2     & 1 \\ 
\hline
~2  &  148 & 112  &  0  &  0  & 36  &  72 & 40  &  2   &    8   &   6   &   0   &   0   &   0    & II,\,1           & 2304      & 3     & 2 \\
\hline
~3  &  130 & 88   &  0  & 12  & 36  &  60 & 22  &  1   &    6   &   6   &   3   &   0   &   0    & III,\,2          & 27648     & 4     & 2 \\
\hline 
~4  &  121 & 76   &  0  & 27  & 27  &  45 & 22  &  1   &    3   &   9   &   0   &   3   &   0    & IV,\,3           & 12288     & 6     & 3 \\
\hline
~5  &  118 & 72   &  8  &  0  & 48  &  56 & 6   &  0   &    8   &   0   &   8   &   0   &   0    &  4               & 20736     & 5     & 2 \\
\hline
~6  &  112 & 64   &  0  & 48  &  0  &  48 & 16  &  1   &    0   &  12   &   0   &   0   &   3    &  V,\,5           & 6912      & 7     & 3 \\
~7  &  112 & 64   &  0  & 0   &  96 &  0  & 16  &  0   &    0   &  16   &   0   &   0   &   0    &  6               & 1728      & 18    & 4 \\
\hline
~8  &  109 & 60   &  4  & 16  &  48 & 36  &  5  &  0   &    4   &  4    &   8   &   0   &   0    &  7               & 248832    & 11    & 3 \\
\hline
~9 &  103 & 52   &  12 & 12  &  48 & 24  &  7  &  0   &    4   &   2   &   6   &   4   &   0    &  8,\,9           & 165888    & 8     & 3 \\
10  &  103 & 52   &  6  & 21  &  45 & 27  &  4  &  0   &    3   &   3   &   9   &   1   &   0    &  10,\,11         & 221184    & 9     & 3 \\
\hline
11 & 100  & 48   &  6  & 26  &  42 & 22  &  4  &  0   &    2   &   4   &   8   &   2   &   0    &  12,\,13         & 995328    & 12    & 3 \\ 
12 & 100  & 48   &  0  & 32  &  48 & 16  &  4  &  0   &    0   &   8   &   8   &   0   &   0    &  14              & 62208     & 24    & 4 \\ 
\hline
13 &  94  & 40   & 27  & 0   &  54 &  0  &  13 &  0   &    4   &   0   &   0   &  12   &   0    &  15              & 6144      & 16    & 4 \\ 
14 &  94  & 40   &  6  & 36  &  36 & 12  &  4  &  0   &    0   &   6   &   6   &   4   &   0    &  17,\,18         & 165888    & 20    & 4 \\ 
15 &  94  & 40   &  8  & 33  &  33 & 19  &  1  &  0   &    1   &   3   &   10  &   1   &   1    &  16,\,20         & 663552    & 10    & 3 \\ 
16 &  94  & 40   &  3  & 36  &  42 & 12  &  1  &  0   &    0   &   4   &   12  &   0   &   0    &  19              & 331776    & 22    & 4 \\ 
\hline
17 &  91  & 36   & 13  & 33  &  27 & 15  &  3  &  0   &    1   &   3   &   6   &   5   &   1    &  21,\,25         & 663552    & 17    & 4 \\ 
18 &  91  & 36   & 15  & 27  &  33 & 13  &  3  &  0   &    1   &   3   &   7   &   3   &   2    &  22,\,24         & 663552    & 13    & 3 \\ 
19 &  91  & 36   &  8  & 36  &  36 &  8  &  3  &  0   &    0   &   4   &   8   &   4   &   0    &  23              & 497664    & 25    & 4 \\ 
\hline
20 &  88  & 32   & 16  & 32  &  24 & 16  &  0  &  0   &    0   &   4   &   8   &   0   &   4    &  26,\,30         & 186624    & $27 \cup 30$ & 4 \\ 
21 &  88  & 32   &  8  & 44  &  24 & 12  &  0  &  0   &    0   &   2   &  10   &   4   &   0    &  28,\,33         & 497664    & 28    & 4 \\ 
22 &  88  & 32   & 10  & 38  &  30 & 10  &  0  &  0   &    0   &   2   &  11   &   2   &   1    &  29,\,32,\,34    & 1990656   & 31    & 4 \\ 
23 &  88  & 32   & 12  & 32  &  36 &  8  &  0  &  0   &    0   &   2   &  12   &   0   &   2    &  27,\,31         & 497664    & 29    & 4 \\ 
24 &  88  & 32   &  8  & 32  &  48 &  0  &  0  &  0   &    0   &   0   &  16   &   0   &   0    &  35              & 124416    & 46    & 4 \\
\hline 
25 &  85  & 28   & 18  & 36  &  24 &  0  &  7  &  0   &    0   &   4   &   0   &  12   &   0    &  36              & 55296     & 19    & 4 \\ 
26 &  85  & 28   & 16  & 36  &  24 &  8  &  1  &  0   &    0   &   2   &   8   &   4   &   2    &  37,\,38         & 1990656   & 23    & 4 \\
27 &  85  & 28   & 12  & 36  &  36 &  0  &  1  &  0   &    0   &   0   &   12  &   4   &   0    &  39              & 55296     & 14    & 3 \\
28 &  85  & 28   & 11  & 40  &  30 &  4  &  0  &  0   &    0   &   0   &   12  &   4   &   0    &  40              & 1658880   & $39 \cup 44$ & 4 \\
\hline
29 &  82  & 24   & 22  & 34  &  18 &  6  &  2  &  0   &    0   &   2   &   5   &   6   &   3    &  41,\,42,\,43    & 995328    & 26 & 4 \\
30 &  82  & 24   & 14  & 48  &  12 &  8  &  0  &  0   &    0   &   0   &   8   &   8   &   0    &  47              & 165888    & 40 & 4 \\
31 &  82  & 24   & 18  & 36  &  24 &  4  &  0  &  0   &    0   &   0   &   10  &   4   &   2    &  45,\,46         & 2985984   & $36 \cup 45$ & 4 \\
32 &  82  & 24   & 22  & 24  &  36 &  0  &  0  &  0   &    0   &   0   &   12  &   0   &   4    &  44              & 165888    & 15 &  3 \\
\hline
33 &  79  & 20   & 22  & 40  &  12 &  4  &  1  &  0   &    0   &   0   &   6   &   8   &   2    &  49,\,50         & 995328    & 35 &  4 \\
34 &  79  & 20   & 25  & 32  &  18 &  4  &  0  &  0   &    0   &   0   &   8   &   4   &   4    &  48              & 1990656   & $37 \cup 42$ & 4  \\
\hline
35 &  76  & 16   & 36  & 24  &  12 &  0  &  4  &  0   &    0   &   2   &   0   &   8   &   6    &  51,\,52         & 82944     & 21 & 4 \\
36 &  76  & 16   & 24  & 48  &  0  &  0  &  4  &  0   &    0   &   0   &   0   &   16  &   0    &  58              & 13824     & 32 & 4 \\
37 &  76  & 16   & 29  & 36  &  6  &  4  &  1  &  0   &    0   &   0   &   4   &    8  &   4    &  56              & 331776    & 34 & 4 \\
38 &  76  & 16   & 32  & 28  &  12 &  4  &  0  &  0   &    0   &   0   &   6   &    4  &   6    &  53,\,55         & 995328    & 43 & 4 \\
39 &  76  & 16   & 32  & 28  &  12 &  4  &  0  &  0   &    0   &   0   &   6   &    4  &   6    &  54,\,57         & 331776    & 38 & 4 \\
\hline
40 &  73  & 12   & 39  & 24  &  6  &  4  &  0  &  0   &    0   &   0   &   4   &    4  &   8    &  59              & 331776    & 41 & 4 \\
\hline
41 &  70  & 8    & 44  & 24  &  0  &  0  &  2  &  0   &    0   &   0   &   0   &    8  &   8    &  60              & 248832    & 48 & 5 \\
\hline
42 &  67  & 4    & 54  & 12  &  0  &  0  &  1  &  0   &    0   &   0   &   0   &    4  &  12    &  61              & 55296     & 33 & 4 \\
\hline
43 &  64  & 0    & 64  & 0   &  0  &  0  &  0  &  0   &    0   &   0   &   0   &    0  &  16    &  62              & 18144     & $47 \cup 49$  & 4,\,5 \\
\hline \hline
\end{tabular}}%
}
\label{tab9}
\end{table}
 
Let us make a brief inspection of Table 9 and highlight the most prominent properties of the hyperplane types. Obviously, a type-one hyperplane is a singular hyperplane and a type-43 hyperplane is an ovoid. Next, let us call a hyperplane of $S_4(3)$ {\it homogeneous} if all its $S_3(3)$'s are of the same type. From Table 9 we discern that homogeneous hyperplanes are of types 7, 24, 36 and 43; there is no homogeneous hyperplane whose $S_3(3)$'s would be all of type $H_1$. At the opposite end of the spectrum, there are hyperplanes of types 15, 17 and 18, which entail all types of $S_3(3)$'s except for that corresponding to $S_{(3)}$'s fully located in it. Interestingly, apart from having the same number of elements, each of them has also one $S_{(3)}$ of type $H_1$ and three $S_{(3)}$'s of type $H_2$. Moreover, all the three feature an odd number of points of order three as, strikingly, there are only two more types (4 and 10) having this property.
Another notable hyperplane types are 22 and 29, where each member stems from {\it three} different types of Veldkamp lines of $S_3(3)$. There are 20 types such that each originates from two different types; each of the remaining 21 types being then the blow-up of a single type of Veldkamp lines of $S_3(3)$. We further observe that there are seven  hyperplane types devoid of points of zeroth order, and the same number of types with a single deep point. It is also worth noting that all those types that feature points of all orders (the first one being type 8, the last one type 37) have some  $S_3(3)$'s of type $H_3$.

\begin{table}[t]
\centering
\caption{The 22 types of Veldkamp points of $S_4(3)$ located on the unique hyperbolic quadric $\mathcal{Q}_0^{+}(15,3)$.  (Compare with Table 6 of \cite{seg2}.)}
\bigskip
\resizebox{\columnwidth}{!}{% 
{\begin{tabular}{|r|r|r|r|r|r|r|r|c|c|c|c|c|c|r|r|c|c|} \hline \hline
%\multicolumn{1}{|c|}{} & \multicolumn{1}{|c|}{} & \multicolumn{1}{|c|}{}  &  \multicolumn{5}{|c|}{}                        & \multicolumn{6}{|c|}{}                      & \multicolumn{1}{|c|}{}  & \multicolumn{1}{|c|}{} & \multicolumn{1}{|c|}{} & \multicolumn{1}{|c|}{} \\
%\cline{4-11}
\multicolumn{1}{|c|}{} & \multicolumn{1}{|c|}{} & \multicolumn{1}{|c|}{}  &  \multicolumn{5}{|c|}{$\#$ of Points of Order} & \multicolumn{6}{|c|}{$\#$ of $S_3(3)$'s of Type} & \multicolumn{1}{|c|}{}  & \multicolumn{1}{|c|}{} & \multicolumn{1}{|c|}{} & \multicolumn{1}{|c|}{} \\
 \cline{4-14}
Tp  &  Pts & Lns  &  0  &  1  &  2  &  3  &  4  &  D   & $H_1$  & $H_2$ & $H_3$ & $H_4$ &  $H_5$ & VL               & Crd       & BS    & W\\
\hline
~1  &  175 & 148  &  0  &  0  &  0  & 108 & 67  &  4   &   12   &   0   &   0   &   0   &   0    & I                & 256       & 2     & 1 \\ 
\hline
~2  &  148 & 112  &  0  &  0  & 36  &  72 & 40  &  2   &    8   &   6   &   0   &   0   &   0    & II,\,1           & 2304      & 3     & 2 \\
\hline
~3  &  130 & 88   &  0  & 12  & 36  &  60 & 22  &  1   &    6   &   6   &   3   &   0   &   0    & III,\,2          & 27648     & 4     & 2 \\
\hline 
~4  &  121 & 76   &  0  & 27  & 27  &  45 & 22  &  1   &    3   &   9   &   0   &   3   &   0    & IV,\,3           & 12288     & 6     & 3 \\
\hline
~6  &  112 & 64   &  0  & 48  &  0  &  48 & 16  &  1   &    0   &  12   &   0   &   0   &   3    &  V,\,5           & 6912      & 7     & 3 \\
~7  &  112 & 64   &  0  & 0   &  96 &  0  & 16  &  0   &    0   &  16   &   0   &   0   &   0    &  6               & 1728      & 18    & 4 \\
\hline
~9 &  103 & 52   &  12 & 12  &  48 & 24  &  7  &  0   &    4   &   2   &   6   &   4   &   0    &  8,\,9           & 165888    & 8     & 3 \\
10  &  103 & 52   &  6  & 21  &  45 & 27  &  4  &  0   &    3   &   3   &   9   &   1   &   0    &  10,\,11         & 221184    & 9     & 3 \\
\hline
13 &  94  & 40   & 27  & 0   &  54 &  0  &  13 &  0   &    4   &   0   &   0   &  12   &   0    &  15              & 6144      & 16    & 4 \\ 
14 &  94  & 40   &  6  & 36  &  36 & 12  &  4  &  0   &    0   &   6   &   6   &   4   &   0    &  17,\,18         & 165888    & 20    & 4 \\ 
15 &  94  & 40   &  8  & 33  &  33 & 19  &  1  &  0   &    1   &   3   &   10  &   1   &   1    &  16,\,20         & 663552    & 10    & 3 \\ 
16 &  94  & 40   &  3  & 36  &  42 & 12  &  1  &  0   &    0   &   4   &   12  &   0   &   0    &  19              & 331776    & 22    & 4 \\ 
\hline 
25 &  85  & 28   & 18  & 36  &  24 &  0  &  7  &  0   &    0   &   4   &   0   &  12   &   0    &  36              & 55296     & 19    & 4 \\ 
26 &  85  & 28   & 16  & 36  &  24 &  8  &  1  &  0   &    0   &   2   &   8   &   4   &   2    &  37,\,38         & 1990656   & 23    & 4 \\
27 &  85  & 28   & 12  & 36  &  36 &  0  &  1  &  0   &    0   &   0   &   12  &   4   &   0    &  39              & 55296     & 14    & 3 \\
28 &  85  & 28   & 11  & 40  &  30 &  4  &  0  &  0   &    0   &   0   &   12  &   4   &   0    &  40              & 1658880   & $39 \cup 44$ & 4 \\
\hline
35 &  76  & 16   & 36  & 24  &  12 &  0  &  4  &  0   &    0   &   2   &   0   &   8   &   6    &  51,\,52         & 82944     & 21 & 4 \\
36 &  76  & 16   & 24  & 48  &  0  &  0  &  4  &  0   &    0   &   0   &   0   &   16  &   0    &  58              & 13824     & 32 & 4 \\
37 &  76  & 16   & 29  & 36  &  6  &  4  &  1  &  0   &    0   &   0   &   4   &    8  &   4    &  56              & 331776    & 34 & 4 \\
38 &  76  & 16   & 32  & 28  &  12 &  4  &  0  &  0   &    0   &   0   &   6   &    4  &   6    &  53,\,55         & 995328    & 43 & 4 \\
39 &  76  & 16   & 32  & 28  &  12 &  4  &  0  &  0   &    0   &   0   &   6   &    4  &   6    &  54,\,57         & 331776    & 38 & 4 \\
\hline
42 &  67  & 4    & 54  & 12  &  0  &  0  &  1  &  0   &    0   &   0   &   0   &    4  &  12    &  61              & 55296     & 33 & 4 \\
\hline \hline
\end{tabular}}%
}
\label{tab10}
\end{table}

To illustrate the degree of similarity with the binary case, in Table 10 there are singled out those 22  types of geometric hyperplanes of $S_4(3)$, featuring $(3^7 + 1)(3^8 - 1)/2$ = 7\,176\,640 members in total, that are in a one-to-one correspondence with the points lying on the unique hyperbolic quadric $\mathcal{Q}_0^{+}(15,3) \subset {\rm PG}(15,3) \subset \mathcal{V}(S_4(3))$ that contains the (image of) $S_4(3)$  and is invariant under its stabilizer group; by comparing Table 10 with Table 5 it can readily be seen that these hyperplanes are blow-ups of those (ordinary) Veldkamp lines of $S_3(3)$ whose cores feature a single point (mod 3).
In addition, seven of them, shown in Table 11, correspond bijectively to the image, furnished by LGr$(4,8)$ (see \cite{hsl}), of the set of $(3+1)(3^2+1)(3^3+1)(3^4+1)$ = 91\,840 generators of the symplectic polar space $\mathcal{W}(7,3) \subset \mathcal{V}(S_3(3))$; note that these are exactly the types whose $S_3(3)$'s are, in complete analogy to the binary case, devoid of both $H_3$'s and $H_5$'s. 

%\bigskip
\begin{table}[pth!]
\centering
\caption{The seven particular types of points of the $\mathcal{Q}_0^{+}(15,3)$ that correspond bijectively to maximal subspaces of the symplectic polar space $\mathcal{W}(7,3)$.  (Compare with Table 7 of \cite{seg2}.)}
\bigskip
\resizebox{\columnwidth}{!}{% 
{\begin{tabular}{|r|r|r|r|r|r|r|r|c|c|c|c|c|c|r|r|c|c|} \hline \hline
%\multicolumn{1}{|c|}{} & \multicolumn{1}{|c|}{} & \multicolumn{1}{|c|}{}  &  \multicolumn{5}{|c|}{}                        & \multicolumn{6}{|c|}{}                      & \multicolumn{1}{|c|}{}  & \multicolumn{1}{|c|}{} & \multicolumn{1}{|c|}{} & \multicolumn{1}{|c|}{} \\
%\cline{4-11}
\multicolumn{1}{|c|}{} & \multicolumn{1}{|c|}{} & \multicolumn{1}{|c|}{}  &  \multicolumn{5}{|c|}{$\#$ of Points of Order} & \multicolumn{6}{|c|}{$\#$ of $S_3(3)$'s of Type} & \multicolumn{1}{|c|}{}  & \multicolumn{1}{|c|}{} & \multicolumn{1}{|c|}{} & \multicolumn{1}{|c|}{} \\
 \cline{4-14}
Tp  &  Pts & Lns  &  0  &  1  &  2  &  3  &  4  &  D   & $H_1$  & $H_2$ & $H_3$ & $H_4$ &  $H_5$ & VL               & Crd       & BS    & W\\
\hline
~1  &  175 & 148  &  0  &  0  &  0  & 108 & 67  &  4   &   12   &   0   &   0   &   0   &   0    & I                & 256       & 2     & 1 \\ 
\hline
~2  &  148 & 112  &  0  &  0  & 36  &  72 & 40  &  2   &    8   &   6   &   0   &   0   &   0    & II,\,1           & 2304      & 3     & 2 \\
\hline
~4  &  121 & 76   &  0  & 27  & 27  &  45 & 22  &  1   &    3   &   9   &   0   &   3   &   0    & IV,\,3           & 12288     & 6     & 3 \\
\hline
~7  &  112 & 64   &  0  & 0   &  96 &  0  & 16  &  0   &    0   &  16   &   0   &   0   &   0    &  6               & 1728      & 18    & 4 \\
\hline
13 &  94  & 40   & 27  & 0   &  54 &  0  &  13 &  0   &    4   &   0   &   0   &  12   &   0    &  15              & 6144      & 16    & 4 \\ 
\hline 
25 &  85  & 28   & 18  & 36  &  24 &  0  &  7  &  0   &    0   &   4   &   0   &  12   &   0    &  36              & 55296     & 19    & 4 \\ 
\hline
36 &  76  & 16   & 24  & 48  &  0  &  0  &  4  &  0   &    0   &   0   &   0   &   16  &   0    &  58              & 13824     & 32 & 4 \\
\hline \hline
\end{tabular}}%
}
\label{tab11}
\end{table}
 
%\newpage
Obviously, as in the preceding case (see Table 4), our projective hyperplane types should be in a one-to-one correspondence with large orbits of the associated $2 \times 2 \times 2 \times 2$  arrays over the three-element Galois field, computed by Bremner and Stavrou \cite{brst}. However, our number of distinct types (43) falls short of theirs (48, disregarding the trivial orbit). Hence, as we also see from Table 9, there are five particular types of hyperplanes such that each should split into two different subtypes. And this indeed happens if we take into account our refined classification of lines of $\mathcal{V}(S_3(3))$ from Sec.\,4.2. For when we check the column `VL' of Table 9 we see that the types of Veldkamp lines that generate geometric hyperplanes corresponding to the union of two BS-orbits are exactly those listed Table 7 at the end of Sec.\,4.2; moreover, as one can readily verify, each subtype has the right cardinality to give the correct number of elements in each particular BS-orbit involved.  By a way of example, let us consider geometric hyperplanes of $S_4(3)$ of type 43, which are generated by lines of $\mathcal{V}(S_3(3))$ of type 62. From Table 7 we see that we have two distinct subtypes, $62a$ and $62b$, having 108 and 648 members, respectively. As each (ordinary) Veldkamp line blows up into $4!=24$ geometric hyperplanes, type 43 splits into subtype $43a$ of cardinality $24 \times 108 = 2592$  and subtype $43b$ having $24 \times 648 = 15\,552$ elements.
And, indeed, from Table 6 of \cite{brst} we see that the large orbit  $\# 47$ has $5184 = 2 \times 2592$ elements and orbit $\# 49$ features $31\,104 = 2 \times 15\,552$ ones. (The difference by a factor of 2 is simply due to the fact that we work with projective spaces, whereas Bremner and Stavrou carried out their computations in the associated vector spaces.)

%however, these splittings are only  in addition to geometry and combinatorics, also group-theoretical methods are employed.

\begin{table}[t]
\centering
\caption{The four prominent types of non-projective points of $\mathcal{V}(S_4(3))$. The Veldkamp line of $S_3(3)$ denoted as V$^{\star}$ consists of the $S_3(3)$ and one of its non-projective ovoids counted three times. (The notation is the same as in Table 9.)}
\bigskip
\resizebox{\columnwidth}{!}{% 
{\begin{tabular}{|r|r|r|r|r|r|r|r|c|c|c|c|c|c|r|r|c|c|} \hline \hline
%\multicolumn{1}{|c|}{} & \multicolumn{1}{|c|}{} & \multicolumn{1}{|c|}{}  &  \multicolumn{5}{|c|}{}                        & \multicolumn{6}{|c|}{}                      & \multicolumn{1}{|c|}{}  & \multicolumn{1}{|c|}{} & \multicolumn{1}{|c|}{} & \multicolumn{1}{|c|}{} \\
%\cline{4-11}
\multicolumn{1}{|c|}{} & \multicolumn{1}{|c|}{} & \multicolumn{1}{|c|}{}  &  \multicolumn{5}{|c|}{$\#$ of Points of Order} & \multicolumn{6}{|c|}{$\#$ of $S_3(3)$'s of Type} & \multicolumn{1}{|c|}{}  & \multicolumn{1}{|c|}{} & \multicolumn{1}{|c|}{} & \multicolumn{1}{|c|}{} \\
 \cline{4-14}
Tp           &  Pts & Lns  &  0  &  1  &  2  &  3  &  4  &  D   & $H_1$  & $H_2$ & $H_3$ & $H_4$ &  $H_5$ & VL                         & Crd       & BS    & W\\
\hline
~6$^{\star}$ &  112 & 64   &  0  & 48  &  0  &  48 & 16  &  1   &    0   &  12   &   0   &   0   &   3    &  V$^{\star}$,\,5$^{\star}$ & 2304      & --    & -- \\
\hline
20$^{\star}$ &  88  & 32   & 16  & 32  &  24 & 16  &  0  &  0   &    0   &   4   &   8   &   0   &   4    &  26$^{\star}$              & 31104       & --    & -- \\ 
\hline
35$^{\star}$ &  76  & 16   & 36  & 24  &  12 &  0  &  4  &  0   &    0   &   2   &   0   &   8   &   6    &  51$^{\star}$              & 13824       & --    & -- \\
\hline
43$^{\star}$ &  64  & 0    & 64  & 0   &  0  &  0  &  0  &  0   &    0   &   0   &   0   &    0  &  16    &  62$^{\star}$              &  7776      & --    & -- \\
\hline \hline
\end{tabular}}%
}
\label{tab12}
\end{table}

\newpage
When it comes to non-projective geometric hyperplanes of $S_4(3)$, those originating from non-projective Veldkamp lines of $S_3(3)$ listed in Table 5 play clearly a distinguished role, because they contain only projective geometric hyperplanes. Their basic properties are given in 
Table 12. Note that none of them features $H_1$ but each is endowed with $H_5$ in some of their 16 $S_3(3)$'s; interestingly, the number of points/lines, as well as that of points of any order, is a multiple of four. All remaining non-projective geometric hyperplanes of $S_4(3)$ will be `generated' by those (full-sized) lines of  $\mathcal{V}(S_3(3))$ that contain non-projective ovoids (see Sec.\,4.2).

Finally, it is worth mentioning that at this point we cannot say much about classification of the lines of $\mathcal{V}(S_4(3))$, even about the projective ones since PG$(15,3)$ itself features  as many as ${21523360 \choose 2}/{4 \choose 2} = 38\,604\,583\,680\,520$ lines, which is well beyond the reach of our currently available computer power.

\section{Conclusion}
We have carried out a detailed, computer-aided analysis of the Veldkamp space of the Segre variety that is the $k$-fold direct product of projective lines of size four, $\mathcal{V}(S_k(3))$, and $k$ runs from 2 to 4. The main result of our work, and the principal departure from the corresponding binary case, is that $\mathcal{V}(S_k(3)) \supset $ PG$(2^k-1,3)$, i.\,e. the existence of non-projective elements. Whereas for $\mathcal{V}(S_2(3))$ such elements are found only among lines (each such line comprising four pairwise-disjoint ovoids), $\mathcal{V}(S_3(3))$ already features both points and lines that do not belong to the associated PG$(7,3)$.
Although non-projective points $\mathcal{V}(S_3(3))$ are found only among ovoids, its non-projective lines exhibit a much greater variety. They form, in fact, two distinct families. A line of the first family  consists of projective points only; there are 2268 such lines, falling into four different types (see Table 5, starred types). A line of the second family contains at least one non-projective point (that is a non-projective ovoid); there are 5400 such lines, forming five distinct types (see Table 6). Curiously enough, one type of the second family, namely $44^{\star}$, is of particular stance here as it serves as a {\it sole} non-projective extension of certain class of lines of the binary $\mathcal{V}(S_3(2))$, which are all projective (see Table 8).
Finally, the four prominent types of non-projective geometric hyperplanes of $\mathcal{V}(S_4(3))$, where  the number of points/lines, as well as that of points of any order, is a multiple of four (see Table 12), certainly deserve further attention as well.

%\newpage
\section*{Acknowledgments}
This work was supported by the Slovak Research and Development Agency under the contract $\#$ SK-FR-2017-0002 and the French Ministry of Europe and Foreign Affairs (MEAE) under the project PHC \v Stef\'anik 2018/40494ZJ. The financial support of both the Slovak VEGA Grant Agency, Project $\#$ 2/0003/16, and the French ``Investissements d'Avenir'' programme, project ISITE-BFC (contract ANR-15-IDEX-03), are gratefully acknowledged as well. We are also indebted to Dr. Petr Pracna for the help with the figures. 

%\vspace*{-.1cm}

%\noindent
%(Last update: 14/5/2018) 


\begin{thebibliography}{10}
\itemsep=-2pt
\bibitem{buec}
F. Buekenhout and A. M. Cohen, Diagram Geometry: Related to Classical Groups and Buildings, Springer, Berlin -- Heidelberg, 2013, Sec. 8.2. 
\bibitem{coh}
A. M. Cohen, Point-line paces related to buildings, in Handbook of Incidence Geometry, F. Buekenhout (ed.), Elsevier, Amsterdam, 1995, 647--737.
\bibitem{near}
M. Saniga, P. L\'evay, M. Planat and P. Pracna, Geometric hyperplanes of the near hexagon L$_3 \times $GQ(2,\,2), Lett. Math. Phys. 91 (2010) 93--104.
\bibitem{seg2}
M. Saniga, H. Havlicek, F. Holweck, M.  Planat,  and P. Pracna, Veldkamp-space aspects of a sequence of nested binary Segre varieties, Annales de l'Institut Henri Poincar\' e D 2 (2015) 309--333;
see also http://arXiv:1403.6714. 
\bibitem{ht}
J. W. P. Hirschfeld and J. A. Thas, General Galois Geometries, Oxford University Press, Oxford, 1991. 
\bibitem{shult}
E. E. Shult, Points and Lines: Characterizing the Classical Geometries, Springer, Berlin -- Heidelberg, 2011. 
\bibitem{bkm}
J. De Beule, A. Klein and K. Metsch, Substructures of finite classical polar spaces, in Current Research Topics in Galois Geometry, J. De Beule and L. Storme (eds.), Nova Science Publishers, New York, 2011, 33--59.
\bibitem{brst}
M.\,R. Bremner and S. Stavrou, Canonical forms of $2 \times 2 \times 2$ and $2 \times 2 \times 2 \times 2$ arrays over $F_2$ and $F_3$, Lin. Multilin. Algebra 61 (2013) 986--997; see also arXiv:1112.0298.
\bibitem{dyck}
W. Dyck, \" Uber Aufstellung und Untersuchung von Gruppe und Irrationalit\"at regul\"arer Riemann'sher Fl\" achen, Math. Annalen 17 (1881) 473--509.
\bibitem{epst}
D. Eppstein, The many faces of the Nauru graph, published online at http://11011110.livejournal.com/124705.html, accessed on October 23, 2016. 
\bibitem{gs}
R. M. Green and M. Saniga, The Veldkamp space of the smallest slim dense near hexagon, Int. J. Geom. Methods Mod. Phys. 10 (2013) 1250082. 
\bibitem{hsl}
F. Holweck, M. Saniga and  P. L\'evay, A notable relation between $N$-qubit and $2^{N-1}$-qubit Pauli groups via binary $LGr (N, 2N)$, Symmetry, Integrability and Geometry: Methods and Applications 10 (2014) 041.
\end{thebibliography}
\end{document}